\documentclass{amsart}
\usepackage{amsmath, amscd, amssymb, amsthm}
\usepackage{bbm}
\usepackage{latexsym}
\usepackage{amsfonts}
\usepackage{graphicx}
\usepackage[all,cmtip]{xy}
\usepackage[colorlinks,linkcolor=blue,breaklinks=blue,urlcolor=blue,citecolor=blue,anchorcolor=blue,pagebackref]{hyperref}%
\usepackage{geometry}
\newtheorem{conjecture}{Conjecture}

\renewcommand*\backref[1]{}
\renewcommand*\backrefalt[4]{ \ifcase #1 \or (cited on page #2) \else (cited on pages #2) \fi}

\newcommand{\be}{\begin{equation}}
\newcommand{\ee}{\end{equation}}
\newcommand{\bea}{\begin{eqnarray}}
\newcommand{\eea}{\end{eqnarray}}

\newcommand{\vs}{\vspace{0.5cm}}

\def\XXint#1#2#3{{\setbox0=\hbox{$#1{#2#3}{\int}$ }
\vcenter{\hbox{$#2#3$ }}\kern-.6\wd0}}

\begin{document}

\title[Canonical metric connections with constant holomorphic sectional curvature]{Canonical metric connections with constant holomorphic sectional curvature}

\author{Shuwen Chen}
\address{Shuwen Chen. School of Mathematical Sciences, Chongqing Normal University, Chongqing 401331, China}
\email{{3153017458@qq.com}}\thanks{Chen is supported by Chongqing graduate student research grant No. CYB240227. The corresponding author Zheng is partially supported by National Natural Science Foundations of China
with the grant nos.\,12141101 and 12471039,  Chongqing Normal University grant 24XLB026, and is supported by the 111 Project D21024.}

\author{Fangyang Zheng}
\address{Fangyang Zheng. School of Mathematical Sciences, Chongqing Normal University, Chongqing 401331, China}
\email{20190045@cqnu.edu.cn; franciszheng@yahoo.com} \thanks{}


\begin{abstract}
We consider the conjecture of Chen and Nie concerning the space forms for canonical metric connections of compact Hermitian manifolds. We verify the conjecture for two special types of Hermitian manifolds: complex nilmanifolds with nilpotent $J$, and non-balanced Bismut torsion-parallel manifolds. 
\end{abstract}

\maketitle

\tableofcontents

\markleft{Shuwen Chen and Fangyang Zheng}
\markright{Canonical metric connections with constant holomorphic sectional curvature}

\section{Introduction and statement of results}\label{intro}

The simplest kind of Riemannian manifolds are the so-called  {\em space forms,} which means complete Riemannian manifolds with constant sectional curvature. Their universal covers are respectively the sphere $S^n$, the Euclidean space ${\mathbb R}^n$, or the hyperbolic space ${\mathbb H}^n$, equipped with (scaling of) the standard metrics.

In the complex case, the sectional curvature of Hermitian manifolds in general can no longer be constant (unless it is flat). Instead one requires the {\em holomorphic sectional curvature} to be constant. When the metric is K\"ahler, one gets the so-called {\em complex space forms,} namely complete K\"ahler manifolds with constant holomorphic sectional curvature. Analogous to the Riemannian case, their universal covers are the complex projective space ${\mathbb C}{\mathbb P}^n$, the complex Euclidean space ${\mathbb C}^n$, or the complex hyperbolic space ${\mathbb C}{\mathbb H}^n$, equipped with (scaling of) the standard metrics.

When a Hermitian metric is not K\"ahler, its curvature tensor does not obey all the K\"ahler symmetries in general. As a result, the holomorphic sectional curvature could no longer determine the entire curvature tensor. So one would naturally wonder about when can the holomorphic sectional curvature be  constant. In this direction, a long-standing conjecture is the following:

\begin{conjecture}[{\bf Constant Holomorphic Sectional Curvature Conjecture}] \label{conj1}
Given any compact Hermitian manifold, if the holomorphic sectional curvature of its Chern (or Levi-Civita) connection is a constant $c$, then when $c\neq 0$ the metric must be K\"ahler (hence a complex space form), while when $c=0$ the metric must be Chern (or Levi-Civita) flat. 
\end{conjecture}
 
Note that when $n\geq 3$, there are  compact Chern flat (or Levi-Civita flat) manifolds that are non-K\"ahler. Compact Chern flat manifolds are compact quotients of complex Lie groups by the classic theorem of Boothby \cite{Boothby}, and compact Levi-Civita flat Hermitian threefolds were classified in \cite{KYZ}, while for dimension $4$ or higher it is still an open question. Note also that the compactness assumption is necessary for the conjecture.

For $n=2$, Conjecture \ref{conj1} was confirmed by Balas and Gauduchon \cite{Balas, BG}, Sato and Sekigawa \cite{SS}, and Apostolov, Davidov, and Mushkarov \cite{ADM} in the 1980s and 1990s. In higher dimensions, the first substantial result towards this conjecture is the result by Davidov, Grantcharov, and Mushkarov \cite{DGM}, in which they showed among other things that the only twistor space with constant holomorphic sectional curvature is the complex space form ${\mathbb C}{\mathbb P}^3$. More recently, Chen, Chen, and Nie \cite{CCN} showed that, for locally conformally K\"ahler manifolds, the conjecture holds provided that the holomorphic sectional curvature is a non-positive constant. The conjecture is also known in some other special cases, see for instance \cite{ChenZ1, LZ, RZ, Tang}. 

Given a Hermitian manifold $(M^n,g)$, besides the Chern connection $\nabla^c$ and Levi-Civita connection $\nabla$, there is another metric connection that is widely studied: the {\em Bismut connection} $\nabla^b$. It is the connection compatible with both the metric $g$ and the almost complex structure $J$, as well as having totally skew-symmetric torsion. Its existence and uniqueness was proved by Bismut in \cite{Bismut}. It was  discovered independently by Strominger \cite{Strominger}, so in some literature it was also called the Strominger connection. 

In our previous work, we tried to extend Conjecture \ref{conj1} to the Bismut connection case, and raised the following conjecture and question \cite[Conjecture 2 and Question 1]{ChenZ}:

\begin{conjecture} \label{conj2}
Given any compact Hermitian manifold, if the holomorphic sectional curvature of its Bismut connection is a non-zero constant, then the metric must be K\"ahler (hence a complex space form). 
\end{conjecture}

\noindent {\bf Question 3.} {\em What kind of compact Hermitian manifolds will have zero Bismut holomorphic sectional curvature but are not Bismut flat?}

\vspace{0.2cm}

The reason for the above splitting is due to the fact that there are examples of compact Hermitian manifolds with vanishing Bismut holomorphic sectional curvature but its Bismut curvature is not identically zero, e.g. the standard Hopf manifolds of dimension $\geq 3$ (or standard Hopf surfaces but with a specially varied metric).  In \cite{ChenZ} we proved Conjecture \ref{conj2} and answered Question 3 for the $n=2$ case, and also answered them in the special case of complex nilmanifolds with nilpotent $J$ and the Bismut K\"ahler-like manifolds (see \cite[Theorems 3 and 4]{ChenZ}).

The three canonical connections on a given Hermitian manifold $(M^n,g)$, namely, $\nabla$ (Levi-Civita), $\nabla^c$ (Chern) and $\nabla^b$ (Bismut), are all equal to each other when $g$ is K\"ahler. While when $g$ is not K\"ahler, they are linearly independent. The line of connections spanned by Chern and Bismut are called {\em Gauduchon connections,} discovered by Gauduchon \cite{Gauduchon1}:
\begin{equation*} \label{eq:Dr}
 D^r = \frac{1+r}{2} \nabla^c + \frac{1-r}{2}\nabla^b, \ \ \ \ \  r\in {\mathbb R}. 
 \end{equation*}
In the literature there are a number of different ways to parametrize these Gauduchon connections. Here we used Gauduchon's approach, so that $D^1=\nabla^c$ is the Chern connection, $D^{-1}=\nabla^b$ is the Bismut connection. Note that in  this parametrization, $D^0=\frac{1}{2}(\nabla^c+\nabla^b)$ is the Hermitian projection of the Levi-Civita connection $\nabla$, often called the {\em Lichnerowicz connection}. Let
\begin{equation*} \label{eq:Drs}
 D^r_s = (1-s) D^r + s \nabla, \ \ \ \ \ \ \ \ \ \ (r,s) \in \Omega := \{ s\neq 1\} \cup \{ (0,1)\} \subset {\mathbb R}^2. 
 \end{equation*}
In the $rs$-plane ${\mathbb R}^2$, the domain $\Omega$ is the cone over the $r$-axis with vertex $(0,1)$, or equivalently, the {\em plane of canonical metric connections} $D^r_s$ are the cone over the line of Gauduchon connections with vertex at the Levi-Civita connection. As is well-known, when the metric $g$ is K\"ahler, all $D^r_s$ coincide, while when $g$ is not K\"ahler,  $D^r_s\neq D^{r'}_{s'}$ for any two distinct points $(r,s)$, $(r',s')$ in $\Omega$.   

In \cite{CN}, H. Chen and X. Nie gave a beautiful characterization of the possible extension of Conjecture \ref{conj1} to the $2$-parameter family of canonical metric connections $D^r_s$. They discovered the particular subset (which will be called the {\em Chen-Nie curve} from now on)
\begin{equation*} \label{eq:Omega0}
\Gamma = \{ (r,s)\in {\mathbb R}^2 \mid (1-r+rs)^2 + s^2 = 4 \} \subset \Omega ,
\end{equation*}
and proved the following theorem (see \cite[Theorem 2.4]{CN}):

\vspace{0.2cm}

\noindent {\bf Theorem 4 (Chen-Nie).}  {\em Let $(M^2,g)$ be a compact Hermitian surface with pointwise constant holomorphic sectional curvature with respect to its $D^r_s$ connection. Then either $g$ must be K\"ahler, or $(r,s)\in \Gamma$ and $(M^2,g)$ is an isosceles Hopf surface equipped with an admissible metric. }

\vspace{0.2cm}

In other words, in order to extend Conjecture \ref{conj1} to the $D^r_s$ connections, one has to exclude the subset $\Gamma$ and address its zero holomorphic sectional curvature case differently just like what we have seen in the Bismut connection case in \cite{ChenZ}. Following their discovery, it is natural to propose the following:

\vspace{0.2cm}

\noindent {\bf Conjecture 5 (Chen-Nie).}  {\em Given any compact Hermitian manifold, if the holomorphic sectional curvature of its $D^r_s$ connection is a non-zero constant, then the metric must be K\"ahler (hence a complex space form). If the holomorphic sectional curvature of its $D^r_s$ connection is zero and $(r,s)\in \Omega \setminus \Gamma$, then $g$ must be $D^r_s$ flat.}

\vspace{0.2cm}

A companion question is the following:

\vspace{0.2cm}

\noindent {\bf Question 6.} {\em For any $(r,s)\in \Gamma$, what kind of compact Hermitian manifolds will have its $D^r_s$ connection being non-flat but  with vanishing holomorphic sectional curvature? }

\vspace{0.2cm}

The aforementioned theorem of Chen and Nie says that Conjecture 5 is true when $n=2$.  Note that the Chen-Nie curve $\Gamma$ lies between the two horizontal lines $s=-2$ and $s=2$, passing through the point $(\frac{1}{3},-2)$, $(-1,2)$, $(-1,0)$, $(3,0)$, and with $r$ approaching $\pm \infty$ when $s\rightarrow 1$. Note that  $D^{-1}_0$ is just the Bismut connection $\nabla^b$, while $D^{3}_0= 2\nabla^c-\nabla^b$ is the reflection point of $\nabla^b$ with respect to $\nabla^c$, sometimes called the {\em anti-Bismut connection}. In the notation of \cite{ZhaoZ},  $D^{\frac{1}{3}}_{-2}=\nabla^-$, $D^{-1}_2=\nabla^+$ are also special connections in the sense that their curvature moves in sync with the curvature of $\nabla^b$, for instance, the curvatures of $\nabla^+$, $\nabla^-$ and $\nabla^b$ will obey all K\"ahler symmetries (namely, be K\"ahler-like) at the same time. Also, $\nabla^++\nabla^c=2\nabla$, that is, $\nabla^+$ is the reflection point of the Chern connection $\nabla^c$ with respect to the Levi-Civita connection $\nabla$. Similarly, $\nabla^-$ is the reflection point of $D^{-\frac{1}{3}}=\frac{1}{3}\nabla^c+\frac{2}{3}\nabla^b$ with respect to $\nabla '=D^0_{-1}$, the so-called {\em anti-Levi-Civita connection}, while  $D^{-\frac{1}{3}}$ is called the {\em minimal connection} since it has the smallest norm of torsion amongst all Gauduchon connections. 

By a beautiful theorem of Lafuente and Stanfield \cite[Theorem A]{LS}, for any Gauduchon connection $D^r$ other than Chern or Bismut, if a compact Hermitian manifold has flat (or more generally, K\"ahler-like) $D^r$, then the metric must be K\"ahler. For $D^r_s$ with $s\neq 0$, it was proved in \cite[Theorem 4]{ZhaoZ} that if $D^r_s$ is other than $\nabla$, $\nabla'$, $\nabla^+$ or $\nabla^-$, and if a Hermitian manifold has flat (or more generally, K\"ahler-like) $D^r_s$, then the metric must be K\"ahler. 

In other words, with the exception of the three canonical connections $\nabla^c$, $\nabla^b$, $\nabla$ and $\nabla'$, $\nabla^+$, $\nabla^-$, for any other $D^r_s$, a compact Hermitian manifold with flat $D^r_s$ must be K\"ahler. 

The main purpose of this short article is to confirm Conjecture 5 for general dimensional Hermitian manifolds in the special case of either complex nilmanifolds (with nilpotent $J$ in the sense of \cite{CFGU}) or non-balanced Bismut torsion-parallel manifolds. Recall that by a {\em complex nilmanifold} here we mean a compact Hermitian manifold $(M^n,g)$ such that its universal cover is a nilpotent Lie group equipped with a left-invariant complex structure and a compatible left-invariant metric. A compact Hermitian manifold $(M^n,g)$ is said to be {\em Bismut torsion-parallel} (or BTP for brevity) if $\nabla^bT^b=0$, where $T^b$ is the torsion of the Bismut connection $\nabla^b$. Also $g$ is said to be {\em balanced} if $d(\omega^{n-1})=0$, where $\omega$ is the K\"ahler form of $g$. Examples and properties of BTP manifolds are discussed in \cite{ZhaoZ22, ZhaoZ24}. Note that any non-K\"ahler Bismut K\"ahler-like (BKL) manifold is always a non-balanced BTP manifold.

\vspace{0.2cm}

\noindent {\bf Theorem 7.}  {\em Let  $(M^n,g)$  be  a compact Hermitian manifold such that for some $(r,s)\in \Omega$, the canonical metric connection $D^r_s$ of $g$ has constant holomorphic sectional curvature $c$.

(1). Assume that $(M^n,g)$ is a complex nilmanifold with nilpotent $J$. If $D^r_s$ is not the Chern connection, then $c=0$, the Lie group is abelian, and $(M^n,g)$ is (a finite undercover of) a flat complex torus. If $D^r_s$ is the Chern connection, then $c=0$, the Lie group is a (nilpotent) complex Lie group, and $(M^n,g)$ is Chern flat.

(2). If $(M^n,g)$ is a non-balanced BTP manifold, then $c=0$ and $(r,s)\in \Gamma$. 

(3). If $(M^3,g)$ is balanced BTP of dimension $3$, then either it is K\"ahler, or it is Chern flat and $c=0$, $(r,s)=(1,0)$ (namely, with $D^r_s$ being the Chern connection), with $M^3$ being a compact quotient of the simple complex Lie group $\mbox{SO}(3, {\mathbb C})$ equipped with the standard metric. }

\vspace{0.2cm}

The above theorem says that Conjecture 5 holds for all complex nilmanifolds with nilpotent $J$, all compact non-balanced BTP manifolds, and all compact balanced BTP threefolds. In particular, since (non-K\"ahler) Bismut K\"ahler-like (BKL) manifolds or Vaisman manifolds are always non-balanced BTP, we know that Conjecture 5 holds for BKL or Vaisman manifolds.

\vspace{0.3cm}

\section{Preliminaries}

First let us recall the definition of {\em sectional curvature} and {\em holomorphic sectional curvature}. Given a  connection $D$ on a differential  manifold $M^n$, its torsion and curvature are respectively defined  by
$$ T^D(x,y)=D_xy-D_yx-[x,y], \ \ \ \ \ R^D_{xy}z=D_xD_yz-D_yD_xz -D_{[x,y]}z, $$
where $x,y,z$ are vector fields on $M$. When $M$ is equipped with a Riemannian metric $g=\langle \cdot , \cdot \rangle $, we could use $g$ to lower the index and turn $R^D$ into a $(4,0)$-tensor (which we still denote by the same letter): 
$$R^D(x,y,z,w)=\langle R^D_{xy}z, \, w\rangle ,$$
where $x,y,z,w$ are vector fields on $M$. Clearly, $R^D$ is skew-symmetric with respect to its first two positions. If $D$ is a {\em metric connection}, namely, $Dg=0$, then $R^D$ is skew-symmetric with respect to its last two positions as well, hence it becomes a bilinear form on $\Lambda^2 T\!M$. Note that the presence of torsion $T^D$ usually will make the bilinear form $R^D$ not symmetric in general. The {\em sectional curvature} of $D$ is defined by
$$  K^D( \pi ) = - \frac{ R^D (x\wedge y, x\wedge y) } { |x\wedge y |^2 } , \ \ \ \ \mbox{where} \ \pi = \mbox{span}_{\mathbb R}\{ x,y\} \subset T_pM.$$ 
Here as usual $|x\wedge y|^2 = |x|^2 |y|^2 - \langle x,y\rangle^2$. It is easy to see that the value $K^D(\pi )$ is independent of the choice of the basis of the $2$-plane $\pi$ in the tangent space $T_pM$. Since the bilinear form $R^D$ may not be symmetric, the values of $K$ in general will not determine the entire $R^D$ (but only the symmetric part of $R^D$). 

Now suppose $(M^n,g)$ is a Hermitian manifold and $D$ is a metric connection. Then besides the sectional curvature $K^D$, one also has the {\em holomorphic sectional curvature} $H^D$, which is the restriction of $K^D$ on those $2$-planes $\pi$ that are $J$-invariant: $J\pi = \pi$. In this case, for any non-zero  $x\in \pi$, $\{ x, Jx\}$ is a basis of $\pi$, so we can rewrite $H^D$ in complex coordinates:
$$ H^D(X) = \frac{R^D(X, \overline{X}, X, \overline{X})}{|X|^4}, \ \ \ \ \ \ \ \ X = x-\sqrt{-1}Jx. $$ 
Let $\{ e_1, \ldots , e_n\}$ be a local frame of type $(1,0)$ complex tangent vector fields on $M^n$, and denote by $R^D_{i\bar{j}k\bar{\ell}}=R^D(e_i, \overline{e}_j , e_k, \overline{e}_{\ell})$ the components of $R^D$ under the frame $e$. We see that
\begin{equation*}
 H^D\equiv c  \ \ \ \Longleftrightarrow \ \ \ \widehat{R}^D_{i\bar{j}k\bar{\ell}} = \frac{c}{2}(\delta_{ij}\delta_{k\ell } + \delta_{i\ell }\delta_{kj}), 
\end{equation*}
where 
$$ \widehat{R}^D_{i\bar{j}k\bar{\ell}}  = \frac{1}{4} \big( R^D_{i\bar{j}k\bar{\ell}}+  R^D_{k\bar{j}i\bar{\ell}} +  R^D_{i\bar{\ell}k\bar{j}} +  R^D_{k\bar{\ell}i\bar{j}} \big) $$
is the symmetrization of $R^D$. 

\vspace{0.1cm}

Next let us recall the structure equations of Hermitian manifolds. Let $(M^n,g)$ be a Hermitian manifold and denote by  $\omega$ the  K\"ahler form associated with $g$. Denote by $\nabla$,  $\nabla^c$, $\nabla^b$ the Levi-Civita, Chern, and  Bismut connection, respectively. Denote by $R$ the curvature of $\nabla$, by $T^c=T$ and $R^c$ the torsion and curvature of $\nabla^c$, and by $T^b$ and $R^b$ the torsion and curvature of $\nabla^b$. Under any local unitary frame $e$, let us write
$$ T^c(e_i,e_k)=\sum_{j=1}^n  T^j_{ik}e_j, \ \ \ \ \ \ 1\leq i,k\leq n. $$
Then $T^j_{ik}$ are the Chern torsion components under $e$. Let $\varphi$ be the coframe of local $(1,0)$-forms dual to $e$, namely, $\varphi_i(e_j)=\delta_{ij}$. Denote by $\theta$, $\Theta$ the matrices of connection and curvature of $\nabla^c$ under $e$. Let  $\tau$ be the  column vector under $e$ of the Chern torsion, namely, $\tau_j=\frac{1}{2} \sum_{i,k}T^j_{ik}\varphi_i\wedge \varphi_k$. Then the structure equations and Bianchi identities are
\begin{equation*}
\label{structure-c} \left\{ \begin{array}{llll} d \varphi = - \ ^t\!\theta \wedge \varphi + \tau,   \\
d  \theta = \theta \wedge \theta + \Theta ,  \\
d \tau = - \ ^t\!\theta \wedge \tau + \ ^t\!\Theta \wedge \varphi, \\
d  \Theta = \theta \wedge \Theta - \Theta \wedge \theta. \end{array} \right.
\end{equation*}
Similarly, denote by $\theta^b$, $\Theta^b$ the matrices of connection and curvature of $\nabla^b$ under $e$. Then
\begin{equation*}
\label{structure-s}
\Theta^b = d\theta^b -\theta^b \wedge \theta^b.
\end{equation*}
Let $\gamma = \nabla^b-\nabla^c$ be the tensor, and for simplicity we will also write $\gamma = \theta^b -\theta$ under $e$. Then by \cite{YZ} we have
\begin{equation} \label{eq:gamma}
 \gamma e_i = \sum_j \gamma_{ij} e_j = \sum_{j,k} \big( T^j_{ik}\varphi_k - \overline{T^i_{jk} } \, \overline{\varphi}_k \big) e_j. 
\end{equation}
Also, following the notation of \cite{YZ,YZ1}, the Levi-Civita connection is given by
\begin{equation*}
 \nabla e_i = \sum_j  \big( ( \theta_{ij} +\frac{1}{2}\gamma_{ij} )   e_j  +  \overline{\beta_{ij}} \,\overline{e}_j \big) , \ \ \ \mbox{where} \ \ \  \beta_{ij} = \frac{1}{2} \sum_k \overline{T^k_{ij}} \varphi_k.
\end{equation*}
Therefore, if we write $e=\,^t\!(e_1, \ldots , e_n)$ and $\varphi = \,^t\!( \varphi_1, \ldots , \varphi_n)$ as column vectors, then under the frame $\,^t\!(e,\overline{e})$ the matrices of connection and curvature of $\nabla$ are given by  
$$ \hat{\theta } = \left[ \begin{array}{ll} \theta_1 & \overline{\beta } \\ \beta & \overline{\theta_1 }  \end{array} \right] , \ \  \ \  \widehat{\Theta } = \left[ \begin{array}{ll} \Theta_1 & \overline{\Theta}_2  \\ \Theta_2 & \overline{\Theta}_1   \end{array} \right], \ \ \ \ \theta_1=\theta + \frac{1}{2}\gamma, \ \ \ \ $$
 where
\begin{equation*}
\label{structure-r}
\left\{ \begin{array}{lll} \Theta_1 = d\theta_1 -\theta_1 \wedge \theta_1 -\overline{\beta} \wedge \beta,  \\
\Theta_2 = d\beta - \beta \wedge \theta_1 - \overline{\theta_1 } \wedge \beta,   \\
d\varphi = - \ ^t\! \theta_1 \wedge \varphi - \ ^t\! \beta
\wedge \overline{\varphi } . \end{array} \right.
\end{equation*}

As is well-known, the entries of the curvature matrix $\Theta$ are all $(1,1)$-forms, while the entries of the column vector $\tau $ are all $(2,0)$-forms, under any frame $e$. Since 
$$ \Theta_{ij} = \sum_{k,\ell =1}^n R^c_{k\overline{\ell} i \overline{j}} \varphi_k \wedge \overline{\varphi}_{\ell}, \ \ \ \ \Theta^b_{ij} = \sum_{k,\ell =1}^n \big( R^b_{k\ell i \overline{j}} \,\varphi_k \wedge \varphi_{\ell}  + R^b_{\overline{k}\overline{\ell} i \overline{j}} \,\overline{\varphi}_k \wedge \overline{\varphi}_{\ell} + R^b_{k\overline{\ell} i \overline{j}} \,\varphi_k \wedge \overline{\varphi}_{\ell} \big) , $$
and
\begin{eqnarray*}
 && (\Theta_1)_{ij} = \sum_{k,\ell =1}^n \big( R_{k\ell i \overline{j}} \,\varphi_k \wedge \varphi_{\ell}  + R_{\overline{k}\overline{\ell} i \overline{j}} \,\overline{\varphi}_k \wedge \overline{\varphi}_{\ell} + R_{k\overline{\ell} i \overline{j}} \,\varphi_k \wedge \overline{\varphi}_{\ell} \big) , \\
 && (\Theta_2)_{ij} = \sum_{k,\ell =1}^n \big( R_{k\ell \overline{i} \overline{j}} \,\varphi_k \wedge \varphi_{\ell}  + R_{k\overline{\ell}\overline{i} \overline{j}} \,\varphi_k \wedge \overline{\varphi}_{\ell} \big).
 \end{eqnarray*}

From the structure equations and Bianchi identities, one gets the relationship between this three curvature tensors (\cite[Lemma 2]{ChenZ}):

\vspace{0.2cm}

\noindent {\bf Lemma 8.}  {\em 
Let $(M^n,g)$ be a Hermitian manifold. Then under any local unitary frame $e$, 
\begin{eqnarray}
R_{k\overline{\ell}i\overline{j}} - R^c_{k\overline{\ell}i\overline{j}} &  = & -\frac{1}{2}T_{ik,\overline{\ell}}^j - \frac{1}{2} \overline{T^i_{j\ell , \overline{k}}}   + \frac{1}{4} \sum_r \big(  T^{r}_{ik}\overline{T^r_{j\ell}} - T^{j}_{kr}\overline{T^i_{\ell r}} - T^{\ell}_{ir}\overline{T^k_{j r}}  \big) \label{RRc}, \nonumber \\
R^b_{k\overline{\ell}i\overline{j}} - R^c_{k\overline{\ell}i\overline{j}} &  = &   -T_{ik,\overline{\ell}}^j -  \overline{T^i_{j\ell , \overline{k}}}   + \sum_r \big(   T^{r}_{ik}\overline{T^r_{j\ell}} - T^{j}_{kr}\overline{T^i_{\ell r}}  \big)  \label{eq:RbRc}
\end{eqnarray}
for any $i$, $j$, $k$, $\ell$, where the indices after commas mean covariant derivatives with respect to $\nabla^c$.}

\vspace{0.2cm}

Note that the discrepancy in the coefficients here and  \cite[Lemma 2]{ChenZ} is due to the fact that our $T^j_{ik}$ is twice of that in \cite{ChenZ}. For our later use, we will also need to express the covariant derivatives of torsion in terms of the Bismut connection. Let us use indices after semicolons to denote the covariant derivatives with respect to the Bismut connection. By formula (\ref{eq:gamma}) for $\gamma$ we have
\begin{eqnarray}
T^j_{ik,\ell} - T^j_{ik;\ell} & = &  \, \sum_r \big( T^j_{rk} T^r_{i\ell} + T^j_{ir} T^r_{k\ell} - T^r_{ik} T^j_{r\ell} \big), \label{eq:Tcb}  \nonumber \\
T^j_{ik,\bar{\ell}} - T^j_{ik;\bar{\ell}} & = & \sum_r \big( - T^j_{rk} \overline{T^i_{r\ell}} - T^j_{ir} \overline{T^k_{r\ell}} + T^r_{ik} \overline{T^r_{j\ell}} \big). \label{eq:Tcbbar}
\end{eqnarray}

Next let us examine the curvature components of the canonical metric connections 
$$D^r_s = (1-s)D^r+s\nabla, \ \ \ D^r=\frac{1+r}{2}\nabla^c + \frac{1-r}{2}\nabla^b, \ \ \ \ (r,s) \in \Omega . $$
For convenience, in the following we will fix an arbitrary point $(r,s)\in \Omega$ and write $D$ for $D^r_s$. Under the local unitary frame $e$, we have
$$ De_i = \sum_j \big( \theta^D_{ij} e_j + s \overline{\beta_{ij}} \,\overline{e}_j \big) , \ \ \ \ \ \theta^D = \theta + \frac{1}{2}(1-r+rs)\gamma =\theta + t\gamma , $$
where we wrote $t=\frac{1}{2}(1-r+rs)$. So under the frame $\,^t\!(e,\overline{e})$ the matrices of connection and curvature of $D$ are given by  
$$ \hat{\theta }^D = \left[ \begin{array}{ll} \theta^D & s\overline{\beta } \\ s\beta & \overline{\theta^D}  \end{array} \right] , \ \  \ \  \ \ \widehat{\Theta }^D = \left[ \begin{array}{ll} \Theta_1^D & \overline{\Theta_2^D}  \\ \Theta_2^D & \overline{\Theta^D_1}   \end{array} \right], \ \ \ \ \ \ \theta^D=\theta + t\gamma, \ \ \ \ $$
 where
\begin{equation*}
\label{structure-r}
\left\{ \begin{array}{ll} \Theta_1^D = d\theta^D -\theta^D \wedge \theta^D - s^2\overline{\beta} \wedge \beta,  \\
\Theta_2^D = s\big( d\beta - \beta \wedge \theta^D - \overline{\theta^D } \wedge \beta \big).  
 \end{array} \right.
\end{equation*}
For any fixed point $p\in M$, by the same proof of \cite[Lemma 4]{YZ}, we may choose our local unitary frame $e$ near $p$ so that $\theta^D|_p=0$. So at the point $p$ we have $\theta|_p = -t\gamma|_p$. Let $\gamma'$ be the $(1,0)$-part of $\gamma$. Then we have $\gamma =\gamma' -\gamma'^{\ast}$ where $\gamma'^{\ast}$ denotes the conjugate transpose of $\gamma'$. Note that we always have $\,^t\!\gamma'\wedge \varphi = -2 \tau$, so by the structure equation we get
$$ \partial \varphi_r =  (\frac{1}{2}-t)\sum_{i,k}T^r_{ik}\varphi_i\wedge \varphi_k,  \ \ \ \ \ \overline{\partial} \varphi_r = t  \sum_{i,k} \overline{T^i_{rk}} \varphi_i \wedge \overline{\varphi}_k,     \ \ \ \ \ \mbox{at the point} \ p. $$
Since $\theta^b|_p= -(t-1)\gamma|_p$, at $p$ we have
\begin{eqnarray}
 \overline{\partial} \gamma'_{k\ell} |_p & = & \overline{\partial}\sum_i T^{\ell}_{ki}\varphi_i  \ = \  \sum_{i} \big( \overline{\partial} (T^{\ell}_{ki}) \wedge \varphi_i + T^{\ell}_{ki} \overline{\partial} \varphi_i \big) \nonumber \\
 & = & \sum_{i,j} \big( - \overline{\partial}_j (T^{\ell}_{ki}) + t \sum_r T^{\ell}_{kr} \overline{T^i_{rj}}    \big) \varphi_i \wedge \overline{\varphi}_j\nonumber \\
 & = & \sum_{i,j} \big(    T^{\ell}_{ik;\bar{j}} +  (t-1)\sum_r \big(  T^{r}_{ik} \overline{T^r_{j\ell}}  - T^{\ell}_{ir} \overline{T^k_{jr} } \big)   -\sum_r T^{\ell}_{kr} \overline{T^i_{jr} }   \big) \varphi_i \wedge \overline{\varphi}_j, \nonumber
\end{eqnarray}
where index after the semicolon stands for covariant derivative with respect to the Bismut connection $\nabla^b$. At $p$, we have
$$ \Theta^D_1 = \Theta +\Phi, \ \ \ \ \ \Phi = td\gamma + t^2 \gamma \wedge \gamma - s^2 \overline{\beta}\wedge \beta. $$
We compute the $(1,1)$-part of $\Phi$ at the point $p$:
\begin{eqnarray}
(\Phi_{k\ell})^{1,1} & = & t\overline{\partial} \gamma'_{k\ell} - t\partial \overline{\gamma'_{\ell k}} -t^2 \sum_r \big( \gamma'_{kr} \wedge \overline{\gamma'_{\ell r}}  + \overline{\gamma'_{r k}}  \wedge \gamma'_{r\ell } \big) -s^2 \sum_r \overline{\beta_{kr}} \wedge \beta_{r\ell } \nonumber\\
& = & \sum_{i,j} \{ tT^{\ell}_{ik;\bar{j}} + t\overline{ T^k_{j\ell ; \bar{i}} } +(t^2-2t)(w-v^{\ell}_i) -t(v^j_i+v^{\ell}_k) -\frac{s^2}{4} v^j_k \} \varphi_i \wedge \overline{\varphi}_j. \nonumber
\end{eqnarray}
Here and from now on we will use these abbreviations:
\begin{equation} \label{eq:vandw}
w=\sum_r T^r_{ik} \overline{T^r_{j\ell }}, \ \  v^j_i= \sum_r T^j_{ir} \overline{T^k_{\ell r}}, \ \  v^{\ell}_i= \sum_r T^{\ell}_{ir} \overline{T^k_{j r}}, \ \  v^j_k= \sum_r T^j_{kr} \overline{T^i_{\ell r}}, \ \  v^{\ell}_k= \sum_r T^{\ell}_{kr} \overline{T^i_{j r}}.  
\end{equation}
We therefore conclude the following:

\vspace{0.2cm}

\noindent {\bf Lemma 9.} {\em 
Let $(M^n,g)$ be a Hermitian manifold. Under any local unitary frame $e$, the curvature of $D=D^r_s$ has components 
$$ R^D_{i\bar{j}k\bar{\ell}} = R^c_{i\bar{j}k\bar{\ell}} + t(T^{\ell}_{ik;\bar{j}} + \overline{ T^k_{j\ell ; \bar{i}} } ) +(t^2-2t)(w-v^{\ell}_i) -t(v^j_i+v^{\ell}_k) -\frac{s^2}{4} v^j_k,$$
for any $1\leq i,j,k,\ell \leq n$, where $t=\frac{1}{2}(1-r+rs)$, $R^c$ is the Chern curvature, $w$ and $v^j_i$ etc. are given by (\ref{eq:vandw}), and indices after the semicolon stand for covariant derivatives with respect to the Bismut connection $\nabla^b$.}
\vspace{0.2cm}

Note that using this short hand notation, by (\ref{eq:Tcbbar}) and (\ref{eq:RbRc}), we get

\vspace{0.2cm}

\noindent {\bf Lemma 10.} {\em 
Given a Hermitian manifold $(M^n,g)$, under any local unitary frame $e$, 
\begin{equation*} \label{eq:RbRc2}
R^b_{i\bar{j}k\bar{\ell}} - R^c_{i\bar{j}k\bar{\ell}} \ =  \ T_{ik;\overline{j}}^{\ell} + \overline{T^k_{j\ell ; \overline{i}}}   + v^{\ell}_i - v^j_i -v^{\ell}_k - w \ \ \  \mbox{for all} \ \  1\leq i,j,k,\ell \leq n,
\end{equation*}
where  $w$ and $v^j_i$ etc. are given by (\ref{eq:vandw}), and the indices after the semicolon stand for covariant derivatives with respect to the Bismut connection $\nabla^b$. }

\vspace{0.2cm}

Since $\widehat{T_{ik;\overline{\ell}}^j} =0$, $\widehat{w}=0$, and 
$ \widehat{v}^j_i=\widehat{v}^j_k=\widehat{v}^{\ell}_i=\widehat{v}^{\ell}_k =\frac{1}{4}(v^j_i + v^j_k +v^{\ell}_i + v^{\ell}_k) :=\widehat{v}$,  we finally end up with the following identity which holds for any Hermitian manifold:
\begin{equation} \label{eq:Rhat}
\widehat{R}^D_{i\bar{j}k\bar{\ell}} \ = \ \widehat{R}^c_{i\bar{j}k\bar{\ell}} -(t^2+\frac{s^2}{4})\widehat{v}\  = \ \widehat{R}^b_{i\bar{j}k\bar{\ell}} + (1-t^2-\frac{s^2}{4})\widehat{v}. 
\end{equation}

In the rest of this section, let us recall a basic formula for {\em Lie-Hermitian manifolds,} which means compact Hermitian manifolds with universal cover $(G,J,g)$, where $G$ is a Lie group equipped with a left-invariant complex structure $J$ and a compatible left-invariant metric $g$. Denote by ${\mathfrak g}$ the Lie algebra of $G$. Then the left-invariant complex structure and metric on $G$ correspond to a complex structure  and a metric  on ${\mathfrak g}$: the former means an almost complex structure $J$ on the vector space ${\mathfrak g}$ satisfying the integrability condition
$$ [x,y] - [Jx,Jy] + J[Jx,y] + J[x,Jy] =0 \ \ \ \ \ \ \forall \ x,y \in {\mathfrak g}, $$
while the latter means an inner product $g=\langle , \rangle $ on ${\mathfrak g}$ such that $\langle Jx,Jy\rangle = \langle x,y\rangle$ for any $x,y\in {\mathfrak g}$. Denote by ${\mathfrak g}_{\mathbb C}$ the complexification of ${\mathfrak g}$, and by ${\mathfrak g}^{1,0}$ its $(1,0)$-part, namely,
$$  {\mathfrak g}^{1,0} = \{ x-\sqrt{-1}Jx \mid x \in {\mathfrak g} \}. $$
By a {\em unitary frame} of ${\mathfrak g}$ we mean a basis $e=\{ e_1, \ldots , e_n\}$ of the complex vector space $  {\mathfrak g}^{1,0}\cong {\mathbb C}^n$ so that $ \langle e_i, \overline{e}_j\rangle = \delta_{ij}$ for any $1\leq i,j\leq n$. Here we assume that $G$ has real dimension $2n$ and we have extended $\langle , \rangle$ bilinearly over ${\mathbb C}$. Following the notation of \cite{VYZ, YZ1, ZZ1} we will let 
\begin{equation} \label{CandD}
[e_i, e_j] = \sum_k C^k_{ij} e_k, \ \ \ \ \ \ \ \ [e_i, \overline{e}_j] = \sum_k \big( \overline{D^i_{kj}} e_k - D^j_{ki} \overline{e}_k \big), \ \ \ \ \ \forall \ 1\leq i,j\leq n. 
\end{equation}
If we denote by $\varphi$ the coframe dual to $e$, then the structure equation takes the form
\begin{equation*}
d\varphi_i = -\sum_{j,k=1}^n \big( \frac{1}{2} C^i_{jk}\varphi_j\wedge \varphi_k  + \overline{D^j_{ik}} \varphi_j \wedge \overline{\varphi}_k \big)
\end{equation*}
and the first Bianchi identity which coincides  with the Jacobi identity becomes
\begin{equation*}
\left\{ \begin{split}  \sum_{r=1}^n \big( C^r_{ij}C^{\ell}_{rk} + C^r_{jk}C^{\ell}_{ri} + C^r_{ki}C^{\ell}_{rj} \big) \ = \ 0,   \hspace{3.2cm}\\
 \sum_{r=1}^n \big( C^r_{ik}D^{\ell}_{jr} + D^r_{ji}D^{\ell}_{rk} - D^r_{jk}D^{\ell}_{ri} \big) \ = \ 0,  \hspace{3cm} \\
 \sum_{r=1}^n \big( C^r_{ik}\overline{D^{r}_{j\ell }} - C^j_{rk}\overline{D^{i}_{r\ell }} + C^j_{ri}\overline{D^{k}_{r\ell }}  - D^{\ell}_{ri}\overline{D^{k}_{jr }} + D^{\ell}_{rk} \overline{ D^{i}_{jr }}  \big) \ = \ 0 ,  \end{split} \right.
\end{equation*}
for any $1\leq i,j,k,\ell \leq n$. Under the frame $e$, the Chern connection form and Chern torsion components are
\begin{equation*}
\theta_{ij}=\sum_{k=1}^n \big( D^j_{ik}\varphi_k  -\overline{D^i_{jk}} \overline{\varphi_k} \big), \ \ \ \     T_{ik}^j = - C_{ik}^j - D_{ik}^j + D_{ki}^j.   \label{T}
\end{equation*}
From this, we get the expression for $\gamma$ and the Bismut connection matrix
\begin{equation*}
\theta^b_{ij}=\sum_{k=1}^n \{ (-C^j_{ik} + D^j_{ki}) \varphi_k  + (\overline{ C^i_{jk}} - \overline{D^i_{kj}} ) \, \overline{\varphi}_k \} .
\end{equation*}
Following \cite{ChenZ}, we can take the $(1,1)$-part in $\Theta^b=d\theta^b - \theta^b  \wedge \theta^b $ and obtain
\begin{eqnarray}
 R^b_{k\overline{\ell}i\overline{j}} & = & \big(C^r_{ik}\overline{C^r_{j\ell } } - C^j_{rk}\overline{C^i_{r\ell } }\big) - \big( C^r_{ik} \overline{D^r_{\ell j} } + \overline{C^r_{j\ell} } D^r_{ki}  \big) + \big( C^j_{ir} \overline{D^k_{r\ell } } - C^j_{kr} \overline{ D^i_{\ell r} } \big)  \nonumber \\
&& + \ \big( \overline{ C^i_{jr}} D^{\ell}_{rk }  - \overline{C^i_{\ell r} } D^{j}_{kr}  \big) - \big( D^j_{ri} \overline{ D^k_{r\ell }} + D^{\ell}_{rk}  \overline{D^i_{rj} }  \big) + \big( D^r_{ki} \overline{ D^r_{\ell j}} - D^{j}_{kr}   \overline{D^i_{\ell r} }  \big) \nonumber
\end{eqnarray}
for any $i,j,k,\ell$. Here $r$ is summed up from $1$ to $n$. For our later proofs, we will also need the following result for the symmetrization of $R^b$, which are (33) and (34) from  \cite{ChenZ}:

\vspace{0.2cm}

\noindent {\bf Lemma 11.} {\em Let $(G,J,g)$ be an even-dimensional Lie group equipped with a left-invariant complex structure and a compatible metric. Let $e$ be a unitary frame of ${\mathfrak g}^{1,0}$ and $C$, $D$ be defined by (\ref{CandD}). Then under $e$ we have
\begin{eqnarray}
 4\widehat{R}^b_{k\overline{k}i\overline{i}} & = & - \,\big( |C^i_{rk}|^2 + |C^k_{ri}|^2 + 2\mbox{Re} \big( C^{i}_{ri}\overline{C^k_{rk } }\big)  \big)  + 2\mbox{Re} \big( \overline{C^i_{ir}} (D^k_{rk}-D^k_{kr})  \nonumber  \\
 && +\, \overline{C^k_{kr}} (D^i_{ri}-D^i_{ir}) + \overline{C^i_{kr}} (D^i_{rk}-D^i_{kr}) + \overline{C^k_{ir}} (D^k_{ri}-D^k_{ir}) \big) \label{Lie1}  \nonumber \\
 && -\,   2\big( |D^i_{rk}|^2 + |D^k_{ri}|^2 + 2\mbox{Re} \big( D^{i}_{ri}\overline{D^k_{rk } } \big) \big) + \big( |D^r_{ki}|^2 + |D^r_{ik}|^2 + 2\mbox{Re} \big( D^{r}_{ik}\overline{D^r_{ki } } \big) \big) \nonumber \\
 && - \, \big( |D^i_{kr}|^2 + |D^k_{ir}|^2 + 2\mbox{Re} \big( D^{i}_{ir}\overline{D^k_{kr } } \big) \big)  ,  \label{Lie1} \\
   \widehat{R}^b_{i\overline{i}i\overline{i}} & = &     - |C^i_{ri}|^2  + 2\mbox{Re} \big( \overline{C^i_{ir}} (D^i_{ri}-D^i_{ir})\big) -2  |D^i_{ri}|^2 + |D^r_{ii}|^2 - |D^i_{ir}|^2 . \label{Lie2}
\end{eqnarray}
  }
  
  \vspace{0.2cm}

Finally let us recall this famous result of Salamon \cite[Theorem 1.3]{Salamon}:

\vspace{0.2cm}

\noindent {\bf Theorem 12 (Salamon).} {\em  Let $G$ be a nilpotent Lie group of dimension $2n$ equipped with a left invariant complex structure. Then there exists a coframe $\varphi =\{ \varphi_1, \ldots , \varphi_n\}$ of left invariant $(1,0)$-forms on $G$ such that $$ d\varphi_1 =0, \ \ \ d\varphi_i = {\mathcal I} \{\varphi_1, \ldots , \varphi_{i-1}\} , \ \ \ \forall \ 2\leq i\leq n, $$
where ${\mathcal I}$ stands for the ideal in  exterior algebra of the complexified cotangent bundle generated by those $(1,0)$-forms.
}

\vspace{0.2cm}

Note that when $g$ is a compatible left-invariant metric, clearly one can choose the above coframe $\varphi$ so that it is also unitary. In terms of the structure constants $C$ and $D$ given by (\ref{CandD}), this means
\begin{equation*}
C^j_{ik}=0  \ \ \ \mbox{unless} \ \ j>i \ \mbox{or} \ j>k; \ \ \ \ \ D^j_{ik}=0  \ \ \ \mbox{unless} \ \ i>j.  \label{Salamon}
\end{equation*}

If the complex structure $J$ is nilpotent in the sense of  Cordero, Fern\'{a}ndez, Gray, and Ugarte \cite{CFGU}, then there exists an invariant unitary coframe $\varphi$ so that
\begin{equation}
C^j_{ik}=0  \ \ \ \mbox{unless} \ \ j>i \ \mbox{and} \ j>k; \ \ \ \ \ D^j_{ik}=0  \ \ \ \mbox{unless} \ \ i>j \ \mbox{and} \ i>k.  \label{CFGU}
\end{equation}

We do not know how to prove Conjecture 5 for all nilmanifolds at the present time, but for those with nilpotent $J$ in the sense of \cite{CFGU}, we will be able to confirm the conjecture with the help of (\ref{CFGU}) above.

\vspace{0.3cm}

\section{The standard Hopf manifolds}

Consider the standard (isosceles) Hopf manifold $(M^n,g)$, $n\geq 2$, where the manifold and the K\"ahler form $\omega$ of the metric are given by
\begin{equation} 
\label{Hopf}
 M^n = ({\mathbb C}^n\setminus \{ 0\})/ \langle \phi \rangle, \ \ \ \ \omega = \sqrt{-1} \frac{ \partial \overline{\partial} |z|^2} {|z|^2} , \ \ \ \ \  \phi (z_1, \ldots , z_n) = (a_1z_1, \ldots , a_nz_n) ,
 \end{equation}
where $a_i$ are constants satisfying $0<|a_1|=\cdots =|a_n| < 1$. Here $(z_1,\ldots , z_n)$ denotes the standard Euclidean coordinate of ${\mathbb C}^n$ and $|z|^2$ is shorthand for $|z_1|^2+\cdots +|z_n|^2$. 

Near any given point $p\in M^n$, $(z_1, \ldots , z_n)$ gives a local holomorphic coordinate system, under which the metric $g$ has components
$ g_{k\bar{\ell}} = \frac{1}{|z|^2} \delta_{k\ell }$. If we let  $e_i=|z|\partial_i$, where $\partial_i = \frac{\partial}{\partial z_i}$, then $e$ becomes a local unitary frame of $(M^n,g)$. The Chern  curvature components are
$$ R^c_{i\bar{j}k\bar{\ell}} = |z|^4 R^c(\partial_i, \overline{\partial}_j, \partial_k, \overline{\partial}_{\ell}) = |z|^4\big( -\partial_i \overline{\partial}_j g_{k\bar{\ell}} + \sum_{r,s} \partial_ig_{k\bar{r}} \,\overline{\partial}_j g_{s\bar{\ell}} \, g^{\bar{r}s} \big) = \delta_{ij}\delta_{k\ell} - \frac{\overline{z}_iz_j}{|z|^2}\delta_{k\ell},
$$
for any $1\leq i,j,k,\ell \leq n$. On the other hand, by the defining equation $\partial \omega^{n-1}=-\eta \wedge \omega^{n-1}$, we know that Gauduchon's torsion $1$-form is  $\eta=(n-1)\partial |z|^2$, and the Chern torsion components under the frame $e$ are
$$ T^j_{ik} = \frac{\overline{z}_k}{|z|}\delta_{ij} - \frac{\overline{z}_i}{|z|}\delta_{kj} , \ \ \ \ \ \ \ \forall \ 1\leq i,j,k\leq n. $$ 
From this, we compute
$$ v^j_i=\sum_r T^j_{ir} \overline{T^k_{\ell r}  } =  \delta_{ij}\delta_{k\ell} + \frac{1}{|z|^2} \big( \overline{z}_iz_{\ell} \delta_{kj} - \overline{z}_iz_j \delta_{k\ell } - \overline{z}_kz_{\ell} \delta_{ij} \big) .$$
Taking its symmetrization, we obtain
\begin{equation*}
4 \hat{v} = 2(\delta_{ij}\delta_{k\ell} + \delta_{i\ell }\delta_{kj} ) - \frac{1}{|z|^2} \big( \overline{z}_iz_j \delta_{k\ell } + \overline{z}_kz_{\ell} \delta_{ij} + \overline{z}_iz_{\ell} \delta_{kj} + \overline{z}_kz_j \delta_{i\ell} \big) = 4\widehat{R}^c.
\end{equation*}
So for any canonical metric connection $D^r_s$ where $(r,s)\in \Omega$, by (\ref{eq:Rhat}) we get
$$ \widehat{R}^D = \widehat{R}^c - (t^2+\frac{s^2}{4})\hat{v} = (1-t^2-\frac{s^2}{4})  \widehat{R}^c.$$
In particular, whenever $ t^2+\frac{s^2}{4}=1$, or equivalently, whenever $(r,s)$ belongs to the Chen-Nie curve $\Gamma$, then one would have $ \widehat{R}^D =0$, that is,  the canonical metric connection $D^r_s$ for the standard Hopf manifold will have vanishing holomorphic sectional curvature. 
When $(r,s)\notin \Gamma$, on the other hand, we have
$$  \widehat{R}^D_{1\bar{1}1\bar{1}} = (1-t^2-\frac{s^2}{4})  \widehat{R}^c_{1\bar{1}1\bar{1}} = (1-t^2-\frac{s^2}{4}) (1-\frac{|z_1|^2}{|z|^2}) . $$
Since $n\geq 2$, the right hand side is not a constant function, thus the holomorphic sectional curvature of $D^r_s$ cannot be a constant.

Next we want to check that, for $(r,s)\in \Omega$, when will the canonical metric connection $D^r_s$ of the standard Hopf manifold be flat? To do this, let us fix any $1\leq i,j,k,\ell\leq n$ and introduce the shorthand notations
$$ b_{ij}=\frac{\overline{z}_iz_j}{|z|^2}\delta_{k\ell}, \ \ \ \ b_{i\ell}=\frac{\overline{z}_iz_{\ell}}{|z|^2}\delta_{kj}, \ \ \ \  b_{kj}=\frac{\overline{z}_kz_j}{|z|^2}\delta_{i\ell}, \ \ \ \  b_{k\ell }=\frac{\overline{z}_kz_{\ell}}{|z|^2}\delta_{ij}.$$ 
We compute
$$ w = \sum_r T^r_{ik}\overline{ T^r_{j\ell} } = b_{ij} + b_{k\ell} - b_{i\ell} - b_{kj}. $$
Similarly,
\begin{eqnarray*}
&& v^j_i \, = \,  \delta_{ij}\delta_{k\ell} -b_{ij} - b_{k\ell} + b_{i\ell} , \ \ \ \ 
v^{\ell}_k \, = \,  \delta_{ij}\delta_{k\ell} -b_{ij} - b_{k\ell} + b_{kj} , \\
&& v^{\ell}_i \, = \,  \delta_{i\ell }\delta_{kj} + b_{ij} - b_{i\ell} - b_{kj} , \ \ \ \ 
v^j_k \, = \,  \delta_{i\ell}\delta_{kj} + b_{k\ell} - b_{i\ell} - b_{kj} . 
\end{eqnarray*}
Also, $R^c_{i\bar{j}k\bar{\ell}} = \delta_{ij}\delta_{k\ell} -b_{ij}$. As is well-known (see for instance the proof of \cite[Lemma 4.6]{ZhaoZ24}), the standard Hopf manifolds given by (\ref{Hopf}) are Bismut torsion-parallel  (or equivalently, Vaisman), so when we plug all of these expressions into the formula in Lemma 9 we end up with
\begin{equation*}
R^D_{i\bar{j}k\bar{\ell}} = (1-2t) \delta_{ij} \delta_{k\ell} + (2t-t^2-\frac{s^2}{4}) \delta_{i\ell} \delta_{kj} +(2t-1)b_{ij} + (t^2-\frac{s^2}{4}) b_{k\ell} + (\frac{s^2}{4}-t) (b_{i\ell} + b_{kj}). 
\end{equation*}
In particular, for any $i\neq k$, we have
\begin{eqnarray*}
R^D_{i\bar{i}k\bar{k}} & = &  (1-2t) + (2t-1) \frac{|z_i|^2}{|z|^2}  + (t^2-\frac{s^2}{4}) \frac{|z_k|^2}{|z|^2} ,  \\
R^D_{i\bar{k}k\bar{i}} & = &  (2t-t^2-\frac{s^2}{4}) + (\frac{s^2}{4}-t) \,\big( \frac{|z_i|^2}{|z|^2}  +  \frac{|z_k|^2}{|z|^2}\big) . 
\end{eqnarray*}
Now assume that $R^D=0$. When $n\geq 3$, each of the coefficients on the right hand sides must be zero, so we get $1-2t=0$, $\frac{s^2}{4}=t^2=t$, which lead to a contradiction. This means that $D^r_s$ can never be flat for any $(r,s)$ when $n\geq 3$. When $n=2$, however, we have $|z_i|^2+|z_k|^2=|z|^2$, so the vanishing of $R^D_{i\bar{i}k\bar{k}} $ and $R^D_{i\bar{k}k\bar{i}}$ in this case only give us
$$ t^2-\frac{s^2}{4}-2t + 1 =0, \ \ \ \ t-t^2=0. $$
From this, we conclude that either $t=1$ and $s=0$, or $t=0$ and $s=\pm 2$. Recall that $t=\frac{1}{2}(1-r+rs)$, so we end up with three solutions:
$$ (r,s)=(-1,0), \, (-1,2), \,  (\frac{1}{3},-2). $$
The corresponding connections are $D^{-1}_0=\nabla^b$, $D^{-1}_2=\nabla^+$, and $D^{\frac{1}{3}}_{-2}=\nabla^-$, namely, the Bismut connection, and the two vertices of the Chen-Nie curve $\Gamma$, which lies between the horizontal lines $s=-2$ and $s=2$. Conversely, it is known that when $n=2$, the isosceles Hopf surface is Bismut flat, hence is also $\nabla^+$ and $\nabla^-$ flat (\cite{ZhaoZ}). 

In summary, we have proved the following:

\vspace{0.2cm}

\noindent {\bf Proposition 13.} {\em 
Let $(M^n,g)$ be a standard (isosceles) Hopf manifold given by (\ref{Hopf}), with $n\geq 2$. Then for any $(r,s)\in \Omega \setminus \Gamma$, the canonical metric connection $D^r_s$ cannot have constant holomorphic sectional curvature. For any $(r,s)\in \Gamma$, $D^r_s$ has vanishing holomorphic sectional curvature, but it is not flat except when $n=2$ and $D^r_s$ is $\nabla^b$, $\nabla^+$ or $\nabla^-$. 
}

\vspace{0.2cm}

Recall that the Chen-Nie curve $\Gamma$ is defined by $1=t^2+\frac{s^2}{4}$ where $2t=1-r+rs$.

\vspace{0.3cm}

\section{Proof of Theorem 7}

In this section we will prove the main result, namely Theorem 7. Let us start with the nilmanifold case.

\begin{proof}[{\bf Proof of Theorem 7 for nilmanifolds.}]
Let $(M^n,g)$ be a complex nilmanifold, namely, a compact Hermitian manifold with universal cover $(G,J,g)$, where $G$ is a nilpotent Lie group, $J$ a left-invariant complex structure on $G$, and $g$ a left-invariant metric on $G$ compatible with $J$. We assume that $J$ is nilpotent in the sense of \cite{CFGU}. Now suppose  that for some $(r,s)\in \Omega$, the holomorphic sectional curvature of the canonical metric connection $D=D^r_s$ is a constant $c$. This means that 
$$ \widehat{R}^D_{i\bar{j}k\bar{\ell}} = \frac{c}{2}(\delta_{ij}\delta_{k\ell} + \delta_{i\ell}\delta_{kj}), \ \ \ \ \forall \ 1\leq i,j,k,\ell\leq n, $$
under any unitary frame $e$. By  (\ref{eq:Rhat}), we have
$$ \widehat{R}^b_{i\bar{j}k\bar{\ell}} = \frac{c}{2}(\delta_{ij}\delta_{k\ell} + \delta_{i\ell}\delta_{kj}) + (t^2+\frac{s^2}{4}-1)\widehat{v}.  $$
Therefore,
\begin{equation} \label{Rbhat2}
\widehat{R}^b_{i\bar{i}k\bar{k}} = \frac{c}{2}(1+\delta_{ik}) + (t^2+\frac{s^2}{4}-1) \cdot \frac{1}{4} \sum_r \{  2\mbox{Re}( T^i_{ir} \overline{ T^k_{kr} } ) + |T^k_{ir}|^2 + |T^i_{kr}|^2  \}, \ \ \ \forall \ 1\leq i,k\leq n.
\end{equation}
Choose $i=k$, we get
$$ \widehat{R}^b_{i\bar{i}i\bar{i}} = c + (t^2+\frac{s^2}{4}-1) \sum_r |T^i_{ir}|^2 = c+ (t^2+\frac{s^2}{4}-1) \sum_{r>i} |D^i_{ri}|^2, $$
where in the last equality we used the fact that $T^j_{ik}=-C^j_{ik}-D^j_{ik} + D^j_{ki}$ and (\ref{CFGU}). Comparing the above identity with (\ref{Lie2}) and utilizing (\ref{CFGU}) again, we end up with
$$ - c = (t^2+\frac{s^2}{4}+1) \sum_{r>i} |D^i_{ri}|^2, \ \ \ \ \ \ \ \ \forall \ 1\leq i\leq n. $$
If we choose $i=n$, then the right hand side is vacuum, so we know that $c=0$. Hence $D^i_{\ast i}=0$. From this, we get $T^i_{i\ast}=0$. So (\ref{Rbhat2}) now takes the form
$$ \widehat{R}^b_{i\bar{i}k\bar{k}} =  (t^2+\frac{s^2}{4}-1)  \frac{1}{4}\sum_r \{   |T^k_{ir}|^2 + |T^i_{kr}|^2  \} =\widehat{R}^b_{k\bar{k}i\bar{i}}.$$
Let us assume that $i<k$. Plugging the above into (\ref{Lie1}) and utilizing (\ref{CFGU}), we get
\begin{eqnarray*}
4\widehat{R}^b_{k\bar{k}i\bar{i}} & = & (t^2+\frac{s^2}{4}-1) \,\{  \sum_{r<k}(|C^k_{ir}|^2 + |D^i_{kr}|^2)+\sum_{r>k}(|D^k_{ri}|^2 +|D^i_{rk}|^2) \} \\
& = & \sum_{r<k} \big( -|C^k_{ri}|^2  +|D^r_{ki}|^2 - |D^i_{kr}|^2 \big) -2 \sum_{r>k} \big( |D^i_{rk}|^2 + |D^k_{ri}|^2 \big) .
\end{eqnarray*}
That is,
\begin{equation*}
 (t^2+\frac{s^2}{4}+1) \sum_{r>k}(|D^k_{ri}|^2 +|D^i_{rk}|^2) + (t^2+\frac{s^2}{4})  \sum_{r<k}(|C^k_{ir}|^2 + |D^i_{kr}|^2) = \sum_{r<k}  |D^r_{ki}|^2 .
\end{equation*}
We already know that $D^j_{\ast j}=0$ for any $j$. In particular $D^{\ast}_{2\ast}=0$ by (\ref{CFGU}), so if we take $k=2$ in the above identity, the right hand side would be zero, thus we conclude that $D^2_{\ast 1}=D^1_{\ast 2}=0$. Hence by (\ref{CFGU}) we have $D^{\ast}_{3\ast }=0$. Take $k=3$ in the above identity, again the right hand side is zero which leads to $D^j_{\ast \ell}=0$ whenever $j,\ell \leq 3$. Thus $D^{\ast}_{4\ast}=0$ by (\ref{CFGU}). Repeating this process, we end up with $D=0$. Then by the above identity again, we get $(t^2+\frac{s^2}{4})C=0$. Note that $(t,s)=(0,0)$ means $(r,s)=(1,0)$ or equivalently $D^r_s=\nabla^c$. So when $D^r_s$ is not the Chern connection, we get $C=0$, hence the Lie group $G$ is abelian, and $g$ is K\"ahler and flat. In this case $(M^n,g)$ is a finite undercover of a flat complex torus. When $(t,s)=(0,0)$, the connection $D^r_s$ is the Chern connection. The vanishing of $D$ means that the Lie group $G$ is a complex Lie group, so $g$ is Chern flat. 

In summary, when $D^r_s$ is not the Chern connection, the constancy of holomorphic sectional curvature for $D^r_s$ would imply that the nilpotent group $G$ must be abelian and $g$ is K\"ahler and flat. When $D^r_s$ is the Chern connection, the constancy of Chern holomorphic sectional curvature would imply that $G$ is a (nilpotent) complex Lie group, and $g$ is Chern flat.  This completes the proof of Theorem 7 for the nilmanifold case. 
\end{proof}

We remark that in the nilmanifold case we do not need to exclude any $(r,s)$ values for the metric connection $D^r_s$. As one can see from the above proof, the technical assumption that $J$ is nilpotent is crucial in the argument, and without which we do not know how to complete the proof. It would be an interesting question to answer though.

Next let us prove  Theorem 7 in the BTP case.

\begin{proof}[{\bf Proof of Theorem 7 for non-balanced BTP manifolds.}]
Let $(M^n,g)$ be a compact, non-balanced BTP manifold. Assume that for some $(r,s)\in \Omega$, the canonical metric connection $D=D^r_s$ has constant holomorphic sectional curvature: $H^D=c$. Then under any local unitary frame $e$, we have $\widehat{R}^D_{i\bar{j}k\bar{\ell}} = \frac{c}{2}(\delta_{ij}\delta_{k\ell} + \delta_{i\ell }\delta_{kj} )$ for any $1\leq i,j,k,\ell \leq n$. By (\ref{eq:Rhat}), we get
\begin{equation} \label{eq:Rbhat2}
\widehat{R}^b_{i\bar{j}k\bar{\ell}} = \frac{c}{2}(\delta_{ij}\delta_{k\ell} + \delta_{i\ell }\delta_{kj} ) + (t^2 +\frac{s^2}{4}-1)\hat{v}, \ \ \ \ \ \ \forall \ 1\leq i,j,k,\ell \leq n,
\end{equation}
where $4\hat{v}= v^j_i + v^{\ell}_k + v^j_k + v^{\ell}_i$. Since $g$ is non-balanced BTP, by Definition 1.6 and Proposition 1.7 of \cite{ZhaoZ24}, we know that locally on $M^n$ there always exist the so-called {\em admissible frames,} which means a local unitary frame $e$ such that the Chern torsion components under $e$ enjoy the property  $T^n_{ij}=0$ and $T^j_{in}=\delta_{ij}a_i$ for any $1\leq i,j\leq n$, where $a_i$ are global constants on $M^n$, also, the Bismut curvature components satisfy $R^b_{n\bar{j}k\bar{\ell }}= R^b_{i\bar{j} n\bar{\ell}}=0$ for any $1\leq i,j,k,\ell\leq n$. Let us take $i=j$ and $k=\ell =n$ in (\ref{eq:Rbhat2}). Then we get
$$ 0 = \frac{c}{2}(1+\delta_{in}) + (t^2 +\frac{s^2}{4}-1)\frac{1}{4}|a_i|^2. $$
For $i=n$, since $a_n=0$ we deduce $c=0$. So the above equality becomes 
$$(t^2 +\frac{s^2}{4}-1)|a_i|^2=0$$ 
for each $i$. Since the metric is assumed to be non-balanced, we have $a_1+\cdots +a_{n-1}=\lambda >0$, therefore those $a_i$ cannot be all zero and we must have $t^2 +\frac{s^2}{4}-1=0$. That is, the parameter $(r,s)$ must belong to the Chen-Nie curve $\Gamma$.  
\end{proof}

It remains to deal with the case of balanced BTP threefolds, which relies on the classification result for such threefolds in \cite{ZhaoZ22}. First we need the following:

\vspace{0.2cm}

\noindent {\bf Lemma 14.} {\em 
Let $(M^n,g)$ be a BTP manifold with its $D^{r}_{s}$ connection having constant holomorphic sectional curvature $c$. Then under any local unitary frame  $e$, 
$$ R^b_{   i\bar{j} k\bar{\ell}   } = \frac{c}{2} ( \delta_{ij}\delta_{k\ell} + \delta_{i\ell} \delta_{kj} )  -\frac{1}{2}w +\frac{1}{4}(t^2+\frac{s^2}{4}-3)(v^j_i+v^{\ell}_k) + \frac{1}{4}(t^2+\frac{s^2}{4}+1)(v^{\ell}_i + v^j_k). $$
}

\vspace{0.2cm}

\begin{proof}
For BTP manifolds, by \cite{ZhaoZ22} we know that the Bismut curvature $R^b$ always satisfies the symmetry condition $R^b_{i\bar{j}k\bar{\ell}} = R^b_{k\bar{\ell}i\bar{j}}$ and
$$ Q_{i\bar{j}k\bar{\ell}} : = \, R^b_{i\bar{j}k\bar{\ell}}- R^b_{k\bar{j}i\bar{\ell}} \, = \, - w - v^j_i -  v^{\ell}_k + v^{\ell}_i + v^j_k ,$$
under any local unitary frame. Therefore for BTP manifolds,
\begin{equation} \label{eq:Rbhat}
 \widehat{ R}^b_{i\bar{j}k\bar{\ell}} = \frac{1}{2}(R^b_{i\bar{j}k\bar{\ell}}+ R^b_{k\bar{j}i\bar{\ell}}) = \frac{1}{2}( 2R^b_{i\bar{j}k\bar{\ell}} - Q_{i\bar{j}k\bar{\ell}}) = R^b_{i\bar{j}k\bar{\ell}} + \frac{1}{2}(w + v^j_i +  v^{\ell}_k - v^{\ell}_i - v^j_k).
\end{equation}
Under our assumption $H^D=c$, we have
$ \widehat{R}^D_{i\bar{j}k\bar{\ell} } = \frac{c}{2} ( \delta_{ij}\delta_{k\ell} + \delta_{i\ell } \delta_{kj} )$. On the other hand, by (\ref{eq:Rhat})  we get
$$ \widehat{R}^D- \widehat{R}^b = \frac{1}{4}(1-t^2-\frac{s^2}{4})(v^j_i +  v^{\ell}_k + v^{\ell}_i + v^j_k).$$
Plugging this into (\ref{eq:Rbhat}), we get the desired expression for $R^b$ stated in the lemma.
\end{proof}

Lemma 14 implies that, for a BTP manifold with $D^r_s$ holomorphic sectional curvature being a constant $c$, its Bismut curvature satisfies
\begin{eqnarray}
R^b_{i\bar{i}k\bar{k}} & = &  \frac{c}{2}(1+\delta_{ik})  - \frac{1}{2} \sum_r |T^r_{ik}|^2  + \frac{1}{2} (t^2+\frac{s^2}{4}-3)\mbox{Re} \sum_r T^i_{ir} \overline{T^k_{kr} } \nonumber \\
& &    + \,\frac{1}{4} (t^2+\frac{s^2}{4}+1)\sum_r (|T^i_{kr}|^2 + |T^k_{ir}|^2 )  ,  \ \ \ \ \ \ \ \ \forall\ 1\leq i,k\leq n,  \label{eq:Rbiikk}
\end{eqnarray}
under any local unitary frame $e$.

Next let us recall the classification result from \cite{ZhaoZ22} for balanced BTP threefolds. Let $(M^3,g)$ be a  balanced, non-K\"ahler, compact BTP Hermitian threefold. By the observation in \cite{ZhouZ} and \cite{ZhaoZ22}, for any given point $p\in M$, there always exists a unitary frame $e$ (which will be called {\em special frames} from now on) in a neighborhood of $p$ such that under $e$ the only possibly non-zero Chern torsion components are $a_i=T^i_{jk}$, where $(ijk)$ is a cyclic permutation of $(123)$. Furthermore, each $a_i$ is a global constant on $M^3$, with $a_1=\cdots =a_r>0$, $a_{r+1}=\cdots =0$, where $r=r_B\in \{ 1,2,3\}$  is the rank of the $B$ tensor, which is the global $2$-tensor on any Hermitian manifold defined under any unitary frame by
$ B_{i\bar{j}} = \sum_{k,\ell } T^j_{k\ell} \overline{   T^i_{k\ell }  }$.
The conclusion in \cite{ZhaoZ22} indicates that any compact balanced (but non-K\"ahler) BTP threefold must be one of the following:

\begin{itemize}
\item $r_B=3$, $(M^3,g)$ is a compact quotient of the complex simple Lie group $SO(3,{\mathbb C})$, in particular it is Chern flat.

\item  $r_B=1$, $(M^3,g)$ is the so-called Wallach threefold, namely, $M^3$ is biholomorphic to the flag variety ${\mathbb P}(T_{{\mathbb P}^2} )$ while $g$ is the K\"ahler-Einstein metric $g_0$ minus the square of the null-correlation section. Scale $g$ by a positive constant if necessary, the Bismut curvature matrix under a special frame $e$ is
\begin{equation} \label{eq:Fanotype}
\Theta^b = \left[ \begin{array}{ccc} \alpha+\beta & 0 & 0\\ 0 & \alpha & \sigma \\ 0 & - \bar{\sigma} & \beta \end{array} \right],  \ \ \ \ \ \ \ \ \ \left\{ \begin{array}{lll}  \alpha = \varphi_{1\bar{1}}+ (1\!-\!b)\varphi_{2\bar{2}} + b \varphi_{3\bar{3}} +p \varphi_{2\bar{3}} + \bar{p} \varphi_{3\bar{2}},
 && \\ \beta = \varphi_{1\bar{1}}+ b\varphi_{2\bar{2}} + (1\!-\!b) \varphi_{3\bar{3}} -p \varphi_{2\bar{3}} - \bar{p} \varphi_{3\bar{2}},  & & \\ \sigma = p\varphi_{2\bar{2}} -p \varphi_{3\bar{3}} +q \varphi_{2\bar{3}} + (1\!+\!b) \varphi_{3\bar{2}}, & &
\end{array} \right.
\end{equation}
where $b$ is a real constant, $p,q$ are complex constants, $\varphi$ is the coframe dual to $e$, and we wrote $\varphi_{i\bar{j}}$ for $\varphi_i \wedge \overline{\varphi}_j$ for simplicity.

\item $r_B=2$, in this case $(M^3,g)$ is said to be of {\em middle type}. Again under appropriate scaling of the metric, the Bismut curvature matrix under $e$ becomes
\begin{equation} \label{eq:middletype}
\Theta^b = \left[ \begin{array}{ccc} d\alpha & d\beta_0 & \\ - d\beta_0 & d\alpha & \\ & & 0 \end{array} \right], \ \ \ \ \ \ \ \ \ \left\{ \begin{array}{ll}  d\alpha = \ x(\varphi_{1\bar{1}}+\varphi_{2\bar{2}}) +iy(\varphi_{2\bar{1}}- \varphi_{1\bar{2}}), & \\ d\beta_0 = -iy (\varphi_{1\bar{1}}+\varphi_{2\bar{2}}) + (x\!-\!2)(\varphi_{2\bar{1}}-\varphi_{1\bar{2}}), & \end{array} \right.
\end{equation}
where $x,y$ are real-valued local smooth functions.
\end{itemize}

With this explicit information on Bismut curvature at hand, we are now ready to finish the proof of Theorems 7 for the balanced BTP threefold case.

\begin{proof}[{\bf Proof of Theorem 7 for balanced BTP threefolds.}]
Let $(M^3,g)$ be a  compact balanced BTP Hermitian threefold. Assume that $g$ is not K\"ahler (otherwise by the constancy of holomorphic sectional curvature we already know that the manifold is a complex space form). Suppose that for some $(r,s)\in \Omega$ the canonical metric connection $D^r_s$ of $g$ has constant holomorphic sectional curvature $c$. 

First let us consider the $r_B=3$ case. In this case $g$ is Chern flat, so by (\ref{eq:Rhat}) we have
\begin{equation}  \label{eq:Brank3}
 \frac{c}{2}(\delta_{ij}\delta_{k\ell} + \delta_{i\ell }\delta_{kj} )= \widehat{R}^D_{i\bar{j}k\bar{\ell}} = -(t^2+\frac{s^2}{4})\hat{v}.
 \end{equation}
Now let $e$ be a special frame. Then the only non-zero Chern torsion components under $e$ are $T^1_{23}=T^2_{31}=T^3_{12}=\lambda >0$. Therefore
$ \hat{v}_{i\bar{i}i\bar{i}} = \sum_r |T^i_{ir}|^2 = 0$ for any $1\leq i\leq 3$, 
while for any $1\leq i\neq k \leq 3$ we have
$$ 4\hat{v}_{i\bar{i}k\bar{k}} = \sum_r \{ |T^i_{kr}|^2 + |T^k_{ir}|^2 + 2\mbox{Re} (T^i_{ir} \overline{T^k_{kr}} )\} =\lambda^2+\lambda^2+0=2\lambda^2.$$
If we let $i=j=k=\ell$ in (\ref{eq:Brank3}), then we get $c=0$. If we let $i=j\neq k=\ell$ in (\ref{eq:Brank3}) instead, then we obtain $c=-(t^2+\frac{s^2}{4})\lambda^2$. Thus $t=s=0$. Since $t=\frac{1}{2}(1-r+rs)$, this means that $r=1$. Hence our connection $D^r_s$ is $D^1_0$, which is the Chern connection $\nabla^c$. 

Next let us consider the $r_B=1$ case. In this case $(M^3,g)$ is the Wallach threefold. Under a special frame $e$, the only non-zero Chern torsion component is $T^1_{23}=\lambda >0$. So as in the previous case $\hat{v}_{i\bar{i}i\bar{i}} = \sum_r |T^i_{ir}|^2 = 0$ for any $1\leq i\leq 3$, and for any $1\leq i< k \leq 3$ we have
$$ 4\hat{v}_{i\bar{i}k\bar{k}} = \sum_r \{ |T^i_{kr}|^2 + |T^k_{ir}|^2 + 2\mbox{Re} (T^i_{ir} \overline{T^k_{kr}} )\} =\delta_{i1} \lambda^2.$$
Therefore by Lemma 14 we have 
$R^b_{i\bar{i}i\bar{i}} = c$ for any $1\leq i\leq 3$ and
\begin{equation} \label{eq:Fano2}
R^b_{i\bar{i}k\bar{k}} = \frac{c}{2} -\frac{1}{2}\delta_{i2}\lambda^2 + \frac{1}{4}(t^2+\frac{s^2}{4}+1)\delta_{i1}\lambda^2, \ \ \ \forall \ 1\leq i<k\leq 3. 
\end{equation}
On the other hand, by the formula (\ref{eq:Fanotype}) for the Bismut curvature matrix under $e$, we have
$$ R^b_{1\bar{1}1\bar{1}}=2, \ \ \ R^b_{2\bar{2}2\bar{2}}=1-b, \ \ \ R^b_{1\bar{1}2\bar{2}}=1, \ \ \ R^b_{2\bar{2}3\bar{3}}=b. $$
Therefore $2=1-b=c$ which implies that $c=2$ and $b=-1$, while by (\ref{eq:Fano2}) we get
$$ 1=\frac{c}{2} + \frac{1}{4}(t^2+\frac{s^2}{4}+1)\lambda^2 , \ \ \  b=\frac{c}{2}-\frac{1}{2}\lambda^2. $$
Note that the first equality in the above line gives a contradiction. So in this Fano case the holomorphic sectional curvature of $D^r_s$ for any $(r,s)\in \Omega$ cannot be a constant. 

Finally let us consider the $r_B=2$ case. In this case, under a special frame $e$ the only non-zero Chern torsion components are $T^1_{23}=T^2_{31}=\lambda >0$. So by (\ref{eq:Rbiikk}) we get  $R^b_{i\bar{i}i\bar{i}} =c$ for each $1\leq i\leq 3$, and, for any $1\leq i<k\leq 3$,
$$ R^b_{i\bar{i}k\bar{k}} = \frac{c}{2} -\frac{1}{2}\lambda^2\delta_{k3} + \frac{1}{4}(t^2+\frac{s^2}{4}+1)\lambda^2 (2-\delta_{k3}). $$
From (\ref{eq:middletype}), we get $R^b_{1\bar{1}1\bar{1}}=x$, $R^b_{3\bar{3}3\bar{3}}=0$, hence  $x=c=0$. Also by (\ref{eq:middletype}) we have $R^b_{1\bar{1}2\bar{2}}=x$, thus the above equality gives us
$$ 0 =x= R^b_{1\bar{1}2\bar{2}} = \frac{c}{2} - 0 + \frac{1}{4}(t^2+\frac{s^2}{4}+1)\lambda^2 (2-0) = \frac{1}{2}(t^2+\frac{s^2}{4}+1)\lambda^2 , $$
which is clearly a contradiction. This shows that balanced BTP threefolds of middle type can never have constant holomorphic sectional curvature for $D^r_s$ for any $(r,s)$. This completes the proof the Theorem 7 for the case of balanced BTP threefolds.
\end{proof}

\vspace{0.1cm}

\noindent\textbf{Acknowledgments.} Zheng would like to thank Haojie Chen, Xiaolan Nie, Kai Tang, Bo Yang, and Quanting Zhao for their interest and helpful discussions.

\vs

\end{document}